\documentclass[11pt]{article}
\usepackage{amsfonts}
\usepackage{amssymb}
\usepackage{latexsym}
\usepackage{amscd}
\usepackage[centertags]{amsmath}
\usepackage{xypic}

\newtheorem{theorem}{Theorem}[section]
\newtheorem{corollary}[theorem]{Corollary}

\newtheorem{proposition}[theorem]{Proposition}
\newtheorem{remark}{Remark}
\newtheorem{example}{Example}
\numberwithin{equation}{section}
\newcommand{\Proof}{ \noindent {\bf Proof.}\ }
\newcommand{\qed}{\hfill $\Box$}

\begin{document}

%%%%%%%%%%%%%%%%

\title{Waiting times and stopping probabilities for patterns in Markov chains }
\author{Min-Zhi Zhao$^{1}$ \quad Dong  Xu$^{2}$ \quad Hui-Zeng Zhang$^{3}$\\
                    % Do not remove
%{\footnotesize Supported  by  NSFC (Grant No. 11371317 ).}\\
{\footnotesize 1. School of Mathematical Sciences, Zhejiang
University, Hangzhou 310027,   China,}\\
 {\footnotesize E-mail:zhaomz@zju.edu.cn. Supported  by  NSFC (Grant No. 11371317 ).}\\
  {\footnotesize 2. School of Mathematical Sciences, Zhejiang
University, Hangzhou 310027,   China,}\\
 {\footnotesize E-mail:xudong\_1236@163.com }\\
  {\footnotesize 3. Department of Mathematics,  Hangzhou Normal University, Hangzhou 310036, China  }\\
 {\footnotesize E-mail:zhanghz789@163.com }}

\date{}%{Received: date / Accepted: date}
%%\MSC2010{11K50; 37F35; 37A50; 37A45;}
\maketitle

\begin{abstract}
  Suppose that $\mathcal C$ is a  finite collection of  patterns. Observe a Markov chain until one of the patterns in
$\mathcal C$ occurs as a run. This time is denoted by $\tau$. In this paper, we aim to give an easy way to calculate the mean waiting
time $E(\tau)$ and the stopping probabilities $P(\tau=\tau_A)$ with $A\in\mathcal C$, where $\tau_A$ is the
waiting time until the pattern $A$ appears as a run.
\end{abstract}

{\verb"Keywords":} Pattern, Markov chain, stopping probability, waiting time

{\verb"AMS 2000 subject classifications":}  60J10, 60J22

\section{Introduction}

Suppose that ${\left\{ {{Z_n}} \right\}_{n \ge 1}}$ is a time homogenous Markov chain with finite state space $\Delta$.
A finite sequence   of elements from
$\Delta$ is called a pattern. %In our paper,
We will use a capital letter to denote a pattern.
Use $\mathcal C$ to denote a finite collection of  patterns.
For example, if $\Delta=\{0,1\}$, then $A=1011$ is a pattern while $\mathcal C=\{101,11\}$ is a  finite collection of  patterns.
For a pattern $A$, use $\tau_A$ to denote  the waiting time until $A$ occurs as a run in the sequence $Z_1,Z_2,\cdots$.
Let
$\tau  = \tau_{\mathcal{C}}=\min \{ \tau_A:A\in\mathcal C\}$
be the waiting time till one of the patterns %in $\mathcal C$
appears.
We are interested in the values $E(\tau)$ and
$P(\tau = \tau_A)$ with $A\in\mathcal C$.

In many applications, such as quality control, hypothesis testing, reliability theory  and scan statistics, the distribution of
$\tau$ is very important. In \cite{naus1965distribution} and \cite{naus2002double}, J.~I.~Naus
used a  window with length $w$  to scan a process until time $T$ and then got a scan statistic. The distribution of this scan statistic
can be transformed into the distribution of $\tau_{\mathcal C}$ with some special collection of patterns. For example, if
$\Delta=\{0,1\}$, $w=4$ and the scan statistic is
$$S_{T}=\mathop {\max }\limits_{1 \le i \le T - 3} ({Z_i} + {Z_{i + 1}} +  {Z_{i + 2}}   + {Z_{i + 3}}),$$ then
$S_{T}$ denotes the maximal number of $1$ appears in a window of length $4$ until time $T$.
In this case, $P({S_T} \geq 2) = P({\tau _{\mathcal C} \leq T})$,~ where $\mathcal C=\{11,101,1001\}$.
Another interesting application is Penney-Ante game which is  developed by Walter Penney (see \cite{nishiyama2010pattern}). It is a game with two players.
Player I chooses a triplet of outcomes namely $A$. Then  payer II chooses a different  triplet namely $B$.
  An unbiased coin  is flipped repeatedly until $A$ or $B$ is observed. If $A$ occurs first, then player I wins the game. Otherwise player II wins.
Clearly, the wining probability for player II is $P(\tau_{\mathcal C}=\tau_B)$, where $\mathcal C=\{A,B\}$.  After player I has selected $A$, the most important thing for player II is to
find   an optimal
strategy, that is he should find a triplet $B$ that maximizes his winning probability. In fact, such an optimal
strategy exists (see \cite{brofos2014markov}).

Thanks to its  importance, the occurrence of patterns has
been studied by many people. When $Z_1,Z_2,\cdots$ are i.i.d., S. R. Li \cite{li1980martingale}, H. U. Gerber and S. R. Li \cite{gerber1981occurrence}
used the Martingale method to study the problem. Later
in 1981, L. J. Guibas and A. M. Odlyzko \cite{guibas1981string} used the combinatorial method to obtain the linear equations of
 $E(\tau)$ and $P(\tau=\tau_A)$. When $\{Z_n\}$ is a Markov chain, in 1990, O. Chrysaphinou and S. Papastavridis \cite{chrysaphinou1991occurrence} used the combinatorial method to obtain the linear equations of $E(\tau)$.
In 2002, J. C. Fu and Y. M. Chang  \cite{fu2002probability}  studied $E(\tau)$ by using  Markov chain embedding method.
Later J. Glaz, M. Kulldorff and etc. \cite{glaz2006gambling}, V. Pozdnyakov  \cite{pozdnyakov2008occurrence} introduced  gambling teams
and used Martingale theory to study $E(\tau)$. In 2014, R. J. Gava and D. Salotti \cite{gava2014stopping}  obtained the system of linear equations of $P(\tau=\tau_A)$ with $A\in\mathcal C$ based on  the results of \cite{glaz2006gambling} and \cite{pozdnyakov2008occurrence}.

When $\{Z_n\}$ is a Markov chain, though  the mean waiting time $E(\tau)$
and the stopping probabilities    $P(\tau=\tau_A)$ were obtained in  \cite{glaz2006gambling}, \cite{pozdnyakov2008occurrence} and \cite{gava2014stopping},  the method is complicated.
Briefly speaking, the method is divided into four steps.
Firstly,  define the sets
$\mathcal D^{'} = \left\{ lA:l \in \Delta, A\in\mathcal C  \right\}$ and $\mathcal C^{'} = \left\{ lm{A}:l,m \in \Delta , A\in\mathcal C\right\}$.
Use $\mathcal D^{''}$ and $\mathcal C^{''}$ to denote  the collection of patterns  excluding from
$\mathcal D'$ and from $\mathcal C'$, respectively, the patterns that cannot occur at time $\tau$.
Set
$K^{'} = |\mathcal C| + |\mathcal D^{''}|$ and  $ M^{'}=|\mathcal C^{''}|$.
Secondly, introduce  the  gambling teams, compute the profit matrix $W$ that has $(K^{'}+M^{'})M'$ elements, and compute the
probability of occurrence of the $i$-th ending scenario  with $ i=1,2,\cdots, K'+M'$.
Thirdly, solve a linear system of $M'$ equations in $M'$ variables and then obtain the mean waiting time $E(\tau)$.
Finally, solve about $M'$ linear systems involving $M'$ equations and  $M'$ variables and then get
the stopping probabilities $P(\tau=\tau_A)$.
%with $A\in\mathcal C$.

In this paper, we aim to find a
more easy and effective  method to calculate
 $E(\tau)$ and     $P(\tau=\tau_A)$.
Inspired by the paper \cite{guibas1981string}, we use the combinative probabilistic analysis and
the Markov property.
The main result of our paper is Theorem \ref{th1.1}. It  extend Theorem 3.3 of \cite{guibas1981string} to Markov  case.
Corollary \ref{c2.1} gives a better way to obtain  $E(\tau)$  and  $P(\tau=\tau_A)$ with $A\in\mathcal C$:
 solving only a single  linear system involving $|\Delta|+|\mathcal C|$ equations and  $|\Delta|+|\mathcal C|$ variables.
 The rest of the paper is organized as follows.
In \S2, the main results and the proofs are given.
In \S3,  some examples are discussed.

\section{Main results}
In our paper, suppose that ${\left\{ {{Z_n}} \right\}_{n \ge 1}}$ is a time   homogenous Markov chain with finite state space $\Delta$,
initial distribution ${\mu_i} = P({Z_1} = i)$ and one-step transition probability ${P_{ij}} = P({Z_{n+1}} = j|{Z_n} = i)$.
We will make the following three assumptions.

(A.1) No pattern in  $\mathcal C$ is  a subpattern of another pattern in  $\mathcal C$.

(A.2) For any
$K=K_{1}K_{2} \cdots K_{m} \in \mathcal C$,   $P_{K_{1}K_{2}} \cdots P_{K_{m-1}K_{m}} > 0$.

(A.3) That $P(\tau < \infty)=1$ and $E(\tau) < \infty$.

For a pattern $K$, let $K_{i}$ denote the $i$-th element of $K$,  $|K|$ denote the length of $K$, that is, $K = K_{1}K_{2}\cdots K_{|K|}$.
Let $X_{K}^{(j)} = I_{ \{j\} } (K_{|K|})$.
%=\left\{
 %                                      \begin{array}{ll}
  %                                       1, & \hbox{$ K_{|K|}=j$;} \\
   %                                      0, & \hbox{$ K_{|K|}\not=j$.}
    %                                   \end{array}
     %                                \right.$$
For patterns $K = {K_1} \cdots {K_s}$ and $T = {T_1} \cdots {T_t}$, let
$\{KT\}$ be a subset of $\{1,2,\cdots, s\wedge t\}$ such that an integer  $ k$ is in $\{KT\}$ if and only if
 ${K_{s - k + 1}} \cdots {K_s} = {T_1} \cdots {T_{k}}$. Note that in  \cite{guibas1981string}, the correlation
of $K$ and $T$, denoted by $KT$, is defined as a string over $\{0,1\}$ with the same length as $K$.
The $k$-th bit (from the right) of $KT$ is $1$ if and only if $k\in\{KT\}$.
For example, if $K=101001$ and $T=10010$, then $KT=001001$ and $\{KT\}=\{1,4\}$.

For any $i \in \Delta$ and any pattern $K$, let
$$P_{i \rightarrow K} = P\left((Z_{2},\cdots ,Z_{|K|+1}) = K|Z_{1}=i\right) = P_{iK_{1}}P_{K_{1}K_{2}} \cdots P_{K_{|K|-1}K_{|K|}}.$$
The pattern of length $0$  is  denoted by $\phi$. Set $P_{i \rightarrow \phi}=1$.
 For  any pattern  $K,T \in \mathcal C$, let
%$${\tilde g_{KT}}(z)  =
%\sum\limits_{r \in \{KT\}} z^r\cdot P_{{T_r} \to {T_{r + 1} \cdots T_{\left| T \right|}}}\left/
% P_{{T_1} \to {T_{2}\cdots T_{\left| T \right|}}}\right. .$$ More specifically,
$${\tilde g_{KT}}(z)=\left\{
                       \begin{array}{ll}
                        \left. \sum\limits_{\scriptstyle r \in \{KT\}\hfill\atop
\scriptstyle1 \le r < \left| T \right|\hfill} z^r \cdot P_{T_r \to T_{r + 1} \cdots T_{\left| T \right|}} \right/
 P_{{T_1} \to {T_{2}\cdots T_{\left| T \right|}}}, & \hbox{$K \ne T$;} \\
                        \left. \left( \sum\limits_{\scriptstyle r \in \{KT\}\hfill\atop
\scriptstyle1 \le r < \left| T \right|\hfill} z^r \cdot P_{T_r \to T_{r + 1} \cdots T_{\left| T \right|}}+ z^{\left| T \right|}\right)
\right/
 P_{{T_1} \to {T_{2}\cdots T_{\left| T \right|}}}, & \hbox{$K=T$.}
                       \end{array} \right.$$
For $i \in \Delta$, $K\in\mathcal C$ and $n\geq1$ , define
$${S_i}(n) = P({Z_n} = i,\tau > n) \text{\,\,and\,\,} S_{K}(n) = P(\tau = \tau_{K} = n).$$
Now, define the corresponding  generating functions
$${F_i}(z) = \sum\limits_{n = 1}^\infty  {{S_i}(n) \cdot {z^{-n}}}  \text{\,\,and\,\,}
{f_K}(z) = \sum\limits_{n = 1}^\infty  {{S_K}(n) \cdot {z^{ - n}}}$$
where $z\ge 1$.
Our main result is the  following Theorem.
\begin{theorem}\label{th1.1}
For any $z \geq 1$, the functions ${F_i}(z)$ and ${f_K}(z)$ with $i \in \Delta$ and $K \in \mathcal C $ satisfy the following
system of linear equations:
\begin{equation}\label{1.1}
\left\{ \begin{array}{l}
\sum\limits_{i \in \Delta} {{F_i}(z) \cdot {P_{ij}}}  = z \cdot {F_j}(z) + z \cdot \sum\limits_{K \in \mathcal C }^{} {{f_K}(z) \cdot X_K^{(j)}}  - {\mu_j},~~ j \in \Delta\\
\sum\limits_{i \in \Delta} {{F_i}(z) \cdot {P_{i  {T_1}}}}  = \sum\limits_{K \in \mathcal C }^{} {{f_K}(z) \cdot {{\tilde g}_{KT}}(z)}  - {\mu_{{T_1}}},
~~ T \in \mathcal C
\end{array} \right.
\end{equation}
\end{theorem}

\Proof
Firstly, for $j\in \Delta$ and $n\ge 1$,
\begin{align*}
\sum_{i\in\Delta}{S_i}(n) \cdot P_{ij}
&= P(\tau > n,Z_{n+1} = j)\\
&= P(\tau > n+1,Z_{n+1} = j) + \sum_{K\in\mathcal C}P(\tau = \tau_{K} = n+1 ,Z_{n+1} = j)\\
&= S_j(n + 1) + \sum_{K\in\mathcal C} S_K
(n + 1) \cdot X_{K}^{(j)}
\end{align*}
Thus we have,
$$
\sum_{n = 1}^\infty\sum_{i\in\Delta}  S_i(n) \cdot z^{ - n} \cdot P_{ij}
= z \cdot \sum_{n = 1}^\infty  S_j(n + 1) \cdot z^{- n - 1}
+ z \cdot \sum_{n = 1}^\infty\sum_{K\in\mathcal C}  S_K(n + 1) \cdot z^{- n - 1}\cdot X_{K}^{(j)}.
$$
Note that
$$S_j(1)+\sum_{K\in\mathcal C} S_K(1)\cdot X_{K}^{(j)} =P(Z_1=j)=\mu_j.$$
It follows that
\begin{equation}\label{2.1}
\sum_{i\in\Delta} F_i(z) \cdot P_{ij} = z \cdot F_j(z) + z \cdot\sum_{K\in\mathcal C} f_K(z) \cdot X_K^{(j)} - \mu_j.
\end{equation}

Secondly, for $T\in\mathcal C$ and $i\in \Delta$, define
$$ S_{i,T}(n) =\left\{
                 \begin{array}{ll}
                   0, & \hbox{$n\le |T|$;} \\
                   P(\tau=\tau_T=n,Z_{n - \left| T \right|}= i), & \hbox{$n\ge |T|+1$.}
                 \end{array}
               \right.$$
Define the corresponding generating function $f_{i,T}(z)$ on $z\ge  1$ as
$${f_{i,T}}(z) = \sum\limits_{n = 1}^\infty  {{S_{i,T}}(n) \cdot {z^{ - n}}}.$$
Clearly, when $n \ge \left| T \right| + 1$,
$S_T(n) =\sum_{i\in\Delta} S_{i,T}(n)$. It implies that
$$\sum\limits_{\left| T \right| + 1}^\infty  {{S_T}(n) \cdot {z^{ - n}}}  = \sum_{i
\in \Delta}\sum_{\left| T \right| + 1}^\infty  S_{i,T}(n) \cdot z^{ - n}.$$
Set $P_{T} = P\left((Z_{1},\cdots,Z_{|T|}\right) = T)=\mu_{T_1}\cdot P_{T_1\to T_2\cdots T_{|T|}}$. Then we have
\begin{equation}\label{2.2}
f_T(z) - {z^{ - \left| T\right|}} \cdot {P_T} =\sum_{i\in\Delta} f_{i,T}(z).
\end{equation}

Thirdly, for $T\in\mathcal C$, $i \in \Delta$ and
$n \geq 1$,
\begin{align}
{S_i}(n) \cdot P_{i \to T}&=P\left(\tau > n,Z_n = i, (Z_{n+1}, \cdots ,Z_{n+|T|} ) = T\right)\nonumber\\
&= \sum_{r = 1}^{|T|} P\left(\tau = n+r ,Z_n = i, (Z_{n+1}, \cdots ,Z_{n+|T|} ) = T\right)\nonumber \\
&= \sum_{1 \le r <|T|}\sum_{K \in \mathcal C} P\left(\tau = \tau_K= n+r ,Z_n = i, (Z_{n+1}, \cdots ,Z_{n+|T|} ) = T\right)\nonumber\\
&\,\,\,\,
+ P(\tau = \tau_{T} =n+|T| ,Z_n = i).\label{2.3}
\end{align}
Obviously,
\begin{equation}\label{2.4}
P(\tau = \tau_{T} =n+|T| ,Z_{n} = i) = S_{i,T}(n + \left| T \right|).
\end{equation}
For $1 \leq r <|T|$ and $K \in \mathcal C$, under the condition that $\tau = \tau_{K}= n+r$, we have $(Z_{n+r-|K|+1}, \cdots ,Z_{n+r} ) = K$.
If in addition $Z_{n} = i$ and
$(Z_{n+1}, \cdots ,Z_{n+|T|} ) = T$, then for the reason that $K$ is not a subpattern of $T$ (except that $K$ may be equal to $T$), we have  $|K| \geq r+1,K_{|K|-r+1} \cdots K_{|K|} = T_{1}\cdots T_{r}$ and $K_{|K|-r} = i$,
that is, $r \in \{KT\}$ and $K_{|K|-r} = i$.
Therefore
\begin{align}
&P\left(\tau = \tau_{K}= n+r ,Z_n = i, (Z_{n+1}, \cdots ,Z_{n+|T|} ) = T\right)\nonumber \\
=&P \left (\tau = \tau_{K}= n+r , (Z_{n+r+1}, \cdots ,Z_{n+|T|} ) = (T_{r+1},\cdots,T_{|T|}) \right )\nonumber\\
  &\cdot I_{\{KT\}}(r)\cdot I_{\{i\}}(K_{|K|-r})\nonumber\\
=&{S_K}(n + r) \cdot {P_{{T_r} \to {T_{r + 1} \cdots T_{\left| T \right|}}}}\cdot I_{\{KT\}}(r)\cdot I_{\{i\}}(K_{|K|-r})\label{2.5}
\end{align}
In view of  (\ref{2.3})--(\ref{2.5}), we obtain that
$$S_i(n) \cdot {P_{i \to T}} = \sum_{K \in \mathcal C } \sum_{\scriptstyle r \in \{KT\}\hfill\atop
\scriptstyle1 \le r < \left| T \right|\hfill} S_K (n + r) \cdot P_{T_r \to {T_{r + 1} \cdots T_{\left| T \right|}}} \cdot I_{\{i\}}(K_{|K|-r}) + S_{i,T}(n + \left| T \right|).$$
Consequently,
$$  \sum_{n = 1}^\infty  {{S_i}(n) \cdot {z^{ - n}}}\cdot {P_{i \to T}}  = \sum_{K \in \mathcal C } \sum_{\scriptstyle r \in \{KT\}\hfill\atop
\scriptstyle1 \le r < \left| T \right|\hfill} z^r \cdot P_{{T_r} \to {T_{r + 1} \cdots T_{\left| T \right|}}} \cdot I_{\{i\}}(K_{|K|-r})  \cdot \sum_{n = 1}^\infty  S_K(n + r) \cdot z^{ - n - r} $$
\begin{equation}\label{2.6}
 + {z^{\left| T \right|}} \cdot \sum_{n = 1}^\infty  {{S_{i,T}}(n + \left| T \right|) \cdot {z^{ - n - \left| T \right|}}}.
\end{equation}
Note that for $r \in \{KT\}$ and $1 \leq r <|T|$, we have $r < |K|$. So
$$\sum\limits_{n = 1}^\infty  S_K(n + r) \cdot z^{ - n - r} = f_{K}(z).$$  Hence we can rewrite (\ref{2.6}) as
\begin{equation}\label{2.7}
 F_i(z)\cdot {P_{i \to T}} = \sum_{K \in \mathcal C }f_K(z) \cdot \sum_{\scriptstyle r \in \{KT\}\hfill\atop
\scriptstyle1 \le r < \left| T \right|\hfill} z^r \cdot P_{{T_r} \to {T_{r + 1} \cdots T_{\left| T \right|}}} \cdot  I_{\{i\}}(K_{|K|-r})    + z^{\left| T \right|} \cdot f_{i,T}(z).
\end{equation}
Summing all $i\in\Delta$ gives
\begin{equation}\label{2.8}
\sum_{i \in \Delta}  F_i(z)\cdot {P_{i \to T}}  = \sum_{K \in \mathcal C } f_K(z) \cdot \sum_{\scriptstyle r \in \{KT\}\hfill\atop
\scriptstyle1 \le r < \left| T \right|\hfill} z^r \cdot P_{{T_r} \to {T_{r + 1} \cdots T_{\left| T \right|}}}
  + {z^{\left| T \right|}} \cdot \sum_{i \in \Delta} {{f_{i,T}}(z)}.
\end{equation}

Finally, combining (\ref{2.2}) with (\ref{2.8}), we conclude that
$$\sum_{i\in \Delta}  F_i(z)\cdot {P_{i \to T}}  = \sum_{K \in \mathcal C } {{f_K}(z) \cdot \sum_{\scriptstyle r \in \{KT\}\hfill\atop
\scriptstyle1 \le r < \left| T \right|\hfill} z^r \cdot P_{{T_r} \to {T_{r + 1} \cdots T_{\left| T \right|}}}
  + {z^{\left| T \right|}} \cdot  f_{T}(z)}-P_{T}.$$
 Dividing by $P_{{T_1} \to {T_{2}\cdots T_{\left| T \right|}}}$  on both sides yields that
\begin{equation}\label{2.9}
\sum_{i \in \Delta} F_i(z)\cdot {P_{i  T_1}}  = \sum_{K \in \mathcal C } f_K(z) \cdot \tilde{g}_{KT}(z)-\mu_{T_1}.
\end{equation}
This, together with (\ref{2.1}),  completes the proof.
\qed

\begin{proposition}\label{p2.1}
 The linear system {\rm(\ref{1.1})} is nonsingular.
\end{proposition}
\Proof
%Similar as the proof of Theorem 1 of \cite{guibas1981string}.
W.l.o.g.,
 %Without loss of generality,
 suppose that
$\Delta = \{1,\cdots,m\}$ and  $\mathcal C = \{A,B,\cdots,T\}$. Let
$$Q(z) = \left (                %左括号
  \begin{array}{cccccccc}   %该矩阵一共3列，每一列都居中放置
    P_{11}-z & P_{21} & \cdots & P_{m1} & -zX_{A}^{(1)} & -zX_{B}^{(1)} & \cdots & -zX_{T}^{(1)}\\
    \cdots \\
    P_{1m} & P_{2m} & \cdots & P_{mm}-z & -zX_{A}^{(m)} & -zX_{B}^{(m)} & \cdots & -zX_{T}^{(m)}\\
    P_{1A_{1}} & P_{2A_{1}} & \cdots & P_{mA_{1}} & -\tilde{g}_{AA}(z) & -\tilde{g}_{BA}(z) & \cdots & -\tilde{g}_{TA}(z)\\
    \cdots\\
    P_{1T_{1}} & P_{2T_{1}} & \cdots & P_{mT_{1}} & -\tilde{g}_{AT}(z) & -\tilde{g}_{BT}(z) & \cdots & -\tilde{g}_{TT}(z)\\
  \end{array}
\right ).   $$
Then
we can rewrite  (\ref{1.1}) as
$$Q(z) \left(F_{1}(z),  \cdots,  F_{m}(z),  f_{A}(z), \cdots, f_{T}(z)\right)^T=\left(-\mu_1,\cdots, -\mu_{m},-\mu_{A_{1}},\cdots,-\mu_{T_{1}}\right)^T. $$
%$$Q(z)\left (
 % \begin{array}{c}
  % F_{1}(z)\\
   %\vdots \\
   %F_{m}(z)\\
   %f_{A}(z)\\
   %\vdots\\
   %f_{T}(z)\\
  %\end{array}
%\right )
%= \left (                %左括号
 % \begin{array}{c}
  %-\mu_{1}\\
  %\vdots\\
  %-\mu_{m}\\
  %-\mu_{A_{1}}\\
  %\vdots\\
  %-\mu_{T_{1}}\\
  %\end{array}
%\right )
%$$
Let $\varphi(z)=|Q(z)| $ be the determinant of $Q(z)$.  It suffices to show that $\varphi(z)$ is a nonzero polynomial. Clearly, at the $i$-th row of $Q(z)$ with $1 \leq i \leq m$,
the highest degree  is $1$ and occurs on the diagonal or after the $m$-th column; while
at the $j$-th row with $j \geq m+1$, the highest degree polynomial occurs only on the diagonal. Therefore in the expansion of $\varphi(z)$, the unique highest degree monomial comes from the product of the diagonal terms.
This, together with the fact  the highest degree monomial of $\tilde{g}_{AA}(z)$ is  $\frac{z^{|A|}}{P_{A_{1} \rightarrow  A_{2}\cdots A_{|A|}}}$, implies that
the unique  highest degree monomial of $\varphi(z)$ is
$$(-1)^{m + |\mathcal C|}\frac{1}{P_{A_{1} \rightarrow  A_{2}\cdots A_{|A|}} P_{B_{1} \rightarrow  B_{2}\cdots B_{|B|}} \cdots P_{T_{1} \rightarrow  T_{2}\cdots T_{|T|}} }z^{m + |A|+ \cdots +|T|}.$$
It shows that $\varphi(z)$ is a nonzero polynomial as desired.
\qed

For $i\in\Delta$ and $T\in\mathcal C$, let
$F_{i} = F_{i}(1)$ and $f_{T} = f_{T}(1)$. Then $F_{i}=E\left(\sum\limits_{n<\tau}I_{\{Z_n=i\}}\right)$ is the mean staying time at $i$ before $\tau$,
and $f_{T} = P(\tau = \tau_{T} < \infty)$ is the probability that the pattern $T$ appears first among all the patterns in $\mathcal C$.
 Thus we have $E(\tau ) = 1 + \sum\limits_{i\in\Delta} {{F_i}} $.
Let ${\tilde g_{KT}} = {\tilde g_{KT}}(1)$.
Substituting $z=1$ into  Theorem \ref{th1.1} gives the following Corollary.
\begin{corollary}\label{c2.1}
The following system of linear equations holds:
\begin{equation}\label{1.2}
\left\{ \begin{array}{l}
\sum\limits_{i \in \Delta} {{F_i} \cdot {P_{ij}}}  = {F_j} + \sum\limits_{K \in \mathcal C }^{} {{f_K} \cdot X_K^{(j)}}  - {\mu_j},~j \in \Delta\\
\sum\limits_{i \in \Delta} {{F_i} \cdot {P_{i  {T_1}}}}  = \sum\limits_{K \in \mathcal C }^{} {{f_K} \cdot {{\tilde g}_{KT}}}  - {\mu_{{T_1}}},~T \in \mathcal C
\end{array} \right.
\end{equation}
\end{corollary}

\begin{remark}
{\rm(1)} For $z\ge 1$, define $$F(z)=1+\sum_{i\in\Delta} F_i(z)=\sum\limits_{n = 0}^\infty  P(\tau>n) \cdot z^{-n}$$
and
$$f(z)=\sum_{K\in\mathcal C} f_K(z)=\sum\limits_{n = 1}^\infty  P(\tau=n) \cdot z^{-n}.$$
If we have solved all $f_K(z)$ with $K\in\mathcal C$, then
we can obtain the generating function $f(z)$. In theory, we can obtain the distribution of $\tau$.
Particularly, we can calculate the moments of $\tau$.

{\rm(2)} Theorem \ref{th1.1} is the generalization of Theorem 3.3 of \cite{guibas1981string}.
Summing  all $j\in \Delta$ in the first part  of (\ref{1.1}), we get
\begin{equation}\label{r2.1}
(z-1)\cdot F(z)+z\cdot\sum\limits_{K \in \mathcal C} {f_{K}(z)} = z.
\end{equation}
In the case that $Z_1,Z_2,\cdots$ are i.i.d and $\mu_j>0$ for all $j$, $P_{ij}=\mu_j$ does not depend on $i$.
Dividing by $\mu_{T_1}$ at the both side of the second part of (\ref{1.1}) gives:
\begin{equation}\label{r2.2}
F(z)=\sum\limits_{K \in \mathcal C }^{} {{f_K}(z) \cdot {{\tilde g}_{KT}}(z)}/\mu_{T_1}.
\end{equation}
If we define
$c_{KT}(z) ={\tilde g}_{KT}(z)/(z\cdot \mu_{T_1})= \sum\limits_{r \in \{KT\}} \frac {z^{r - 1}}{\mu_{T_{1}} \cdots \mu_{T_{r}}}$,
then combining (\ref{r2.1}) with (\ref{r2.2}) yields Theorem 3.3 of \cite{guibas1981string}.
Note that the definition of ${c_{KT}}(z)$ in \cite{guibas1981string} has a typo and we  correct it here.

{\rm(3)}  To obtain  $E(\tau)$  and  $P(\tau=\tau_A)$ with $A\in\mathcal C$,  we  only need to
 solve  one   linear system involving $|\Delta|+|\mathcal C|$ equations and  $|\Delta|+|\mathcal C|$ variables.
Compared with the results in  \cite{glaz2006gambling}, \cite{pozdnyakov2008occurrence} and \cite{gava2014stopping}, it is a much
easy and effective way.
\end{remark}

When $|T|=1$  and $T$ is  not a subpattern of $K$,  we must have
$$
\tilde{g}_{KT}(z) =
\left\{ \begin{array}{l}
0,~K \neq T;\\
z,~K = T.
\end{array} \right.
$$
If $j\in \mathcal C$, then  $F_{j}(z) = 0$.
By the above discussion, Theorem \ref{th1.1} yields  the following Corollary.

\begin{corollary}\label{c2.2}
If the lengths of all patterns in $\mathcal C$ are $1$, then the following linear system  holds:
$$
\left\{ \begin{array}{l}
\sum \limits_{ i \notin \mathcal C} {F_{i}(z)\cdot P_{ij}} = z\cdot f_{j}(z) - \mu_{j},~j \in \mathcal C;\\
\sum \limits_{ i \notin \mathcal C} {F_{i}(z)\cdot P_{ij}} = z\cdot F_{j}(z) - \mu_{j},~j \notin \mathcal C.
\end{array} \right.
$$
\end{corollary}

When all pattern contains only one element, we  only need to
 solve  a  linear system involving $|\Delta|$ equations.

\begin{corollary}\label{c2.3}
 Suppose that the first elements of all patterns in $\mathcal C$ are equal and $A$ is any pattern in $\mathcal C$. Then the following linear system   holds:
\begin{equation}\label{eq1}
\left\{ \begin{array}{l}
\sum\limits_{K \in \mathcal C } {{f_K}}  = 1,\\
\sum\limits_{K \in \mathcal C }^{} {{f_K} \cdot ({{\tilde g}_{KT}} - {{\tilde g}_{KA}})}  = 0, ~T \in \mathcal C,~T\not=A.
\end{array} \right.
\end{equation}
\end{corollary}
\Proof
Set $h=A_{1}$. Then $T_{1} = h$ for all $T \in \mathcal C$. In this case, the second part of (\ref{1.2}) can be rewritten as following:
$$
\sum\limits_{i \in \Delta} {{F_i} \cdot {P_{ih}}}  = \sum\limits_{K \in \mathcal C } {{f_K} \cdot {{\tilde g}_{KT}}}  - {\mu_{h}},~T \in \mathcal C.
$$
It shows that for all $T\in\mathcal C$, the values $\sum\limits_{K \in \mathcal C } {{f_K} \cdot {{\tilde g}_{KT}}}$ are the same.
Particularly, $$\sum\limits_{K \in \mathcal C } {{f_K} \cdot {{\tilde g}_{KT}}}=\sum\limits_{K \in \mathcal C } {{f_K} \cdot {{\tilde g}_{KA}}}.$$
This, combining with the fact that $\sum_{K \in \mathcal C } {{f_K}}  = 1$ yields our result.
\qed

When  the first elements of all patterns are equal, namely $h$, the calculation become  more simplified. To solve $f_{K}$ with $K \in \mathcal C$,
it is enough to solve a linear system of
$|\mathcal C|$  equations.
In this case,  the stopping probabilities
 are only related to  the transition probability among those states in $\Delta_1$,
 but neither the initial distribution nor the transition probability $P_{ij}$ with $i$ or $j$ outside $\Delta_1$, where
$\Delta_1$ is the set of elements of patterns in $\mathcal C$.
This is actually true. Intuitively,  all patterns do not occur before the first visiting $h$.
In addition, if the process stays outside $\Delta_1$ and no pattern has occurred, then
the behavior before his next visiting $h$ will not affect the stopping probabilities.

  Sometimes we are interested in when will the distribution of $Z_{\tau}$ is  the same as the initial distribution.
The Corollary below    gives the answer.

\begin{corollary}\label{c2.4}
 Assume that $\{Z_n\}$ is irreducible and has the unique stationary distribution $\pi$.

{\rm (1)} The distribution of $Z_{\tau}$ is  the same as the initial distribution if and only if
there is a constant $c$ such that $F_i=c\cdot\pi_i$ for all $i\in\Delta$.
Actually, $E(\tau)=1+c$ and  $c=({\sum_{K\in\mathcal C}f_K\cdot{\tilde g}_{KT}-\mu_{T_1}})/\pi_{T_1}$
with any given $T\in \mathcal C$.

%{\rm(2)} If the last elements of all patterns are all equal to $t$, then  the distribution of $Z_{\tau}$ is  the same as the %initial distribution if and only if
%$\mu_t=1$.

{\rm(2)}  If the distribution of $Z_{\tau}$ is  the same as the initial distribution, then the following linear system holds:
 \begin{equation}\label{eq2}
\left\{ \begin{array}{l}
\sum\limits_{K \in \mathcal C } {{f_K}}  = 1,\\
\sum\limits_{K \in \mathcal C }^{} {{f_K} \cdot ({{\tilde g}_{KT}} - X_K^{(T_1)})}  = c\cdot \pi_{T_1}, ~T \in \mathcal C
\end{array} \right.
\end{equation}
\end{corollary}

\Proof
By (1) and  Corollary \ref{c2.1},   (2) follows immediately. Thus we only need to prove (1).
 The first part of (\ref{1.2}) shows that  the distribution of $Z_{\tau}$ is  the same as the initial distribution if and only if
\begin{equation}\label{ceq2.1}
\sum_{i\in\Delta} F_i\cdot P_{ij}=F_j,~j\in\Delta.
\end{equation}
Equivalently, there is a constant $c$ such that $F_i=c\cdot \pi_i$ for all $i\in \Delta$.
In this case, $E(\tau)=1+\sum_{i\in\Delta} F_i=1+c$.
By (\ref{ceq2.1}) and the second part of  (\ref{1.2}), we have
$$F_{T_1}=\sum\limits_{K \in \mathcal C }^{} {{f_K} \cdot {{\tilde g}_{KT}}}  - {\mu_{{T_1}}}.$$
It follows that $c=(\sum_{K \in \mathcal C }^{} {{f_K} \cdot {{\tilde g}_{KT}}}  - {\mu_{{T_1}}})/\pi_{T_1}$ as desired.
\qed

\section{Examples}

We begin with the analysis of  Example 1 of \cite{pozdnyakov2008occurrence}.
The mean waiting time and the generating function of $\tau$ are calculated in Example 1 and Example 3 of  \cite{pozdnyakov2008occurrence} respectively,
while the stopping probability is obtained in Example 3.1 of  \cite{gava2014stopping}. We now recalculate all these values  by applying our results.

\begin{example}\label{exmp1}
 Suppose that $\Delta  = \left\{ {1,2,3} \right\},{\rm{\mathcal C = \{ 323,313,33\} }},{\mu_1} = {\mu_2} = {\mu_3} = {1}/{3}$ and the one-step transition probability matrix is
$$P=\left({\begin{array}{*{20}{c}}
{{3 \mathord{\left/
 {\vphantom {3 4}} \right.
 \kern-\nulldelimiterspace} 4}}&0&{{1 \mathord{\left/
 {\vphantom {1 4}} \right.
 \kern-\nulldelimiterspace} 4}}\\
0&{{3 \mathord{\left/
 {\vphantom {3 4}} \right.
 \kern-\nulldelimiterspace} 4}}&{{1 \mathord{\left/
 {\vphantom {1 4}} \right.
 \kern-\nulldelimiterspace} 4}}\\
{{1 \mathord{\left/
 {\vphantom {1 4}} \right.
 \kern-\nulldelimiterspace} 4}}&{{1 \mathord{\left/
 {\vphantom {1 4}} \right.
 \kern-\nulldelimiterspace} 4}}&{{1 \mathord{\left/
 {\vphantom {1 2}} \right.
 \kern-\nulldelimiterspace} 2}}
\end{array}} \right).$$
Let $A = 323$,$B = 313$ and $C = 33$. By calculation, we get
$${\tilde g_{AA}}(z) = z+16z^3,~~{\tilde g_{BA}} = z,~~{\tilde g_{CA}} = z,$$
$${\tilde g_{AB}}(z) = z,~~{\tilde g_{BB}}(z) =z+16z^3,~~{\tilde g_{CB}} = z,$$
$${\tilde g_{AC}} = z,~~{\tilde g_{BC}} = z,~~{\tilde g_{CC}} = z+2z^2.$$
Put these values into  (\ref{1.1}), we get
$$\left(
  \begin{array}{cccccc}
    \frac 34-z & 0 & \frac 14 & 0 & 0 & 0 \\
    0& \frac 34-z & \frac 14 & 0 & 0 & 0 \\
    \frac 14 & \frac 14 & \frac 12-z & -z & -z & -z \\
    \frac 14 & \frac 14 & \frac 12 & -z-16z^3 & -z & -z \\
    \frac 14 & \frac 14 & \frac 12& -z & -z-16z^3 & -z \\
    \frac 14 & \frac 14 & \frac 12 & -z & -z & -z-2z^2 \\
  \end{array}
\right)\left(
         \begin{array}{c}
           F_1(z) \\
           F_2(z) \\
           F_3(z) \\
           f_A(z) \\
           f_B(z)\\
           f_C(z) \\
         \end{array}
       \right)=
\left(
  \begin{array}{c}
    -\frac 13 \\
    -\frac 13\\
    -\frac 13\\
    -\frac 13 \\
    -\frac 13\\
    -\frac 13 \\
  \end{array}
\right).
$$
It is easily seen that
$$f_A(z)=f_B(z)=\frac {F_3(z)}{16z^2}, f_C(z)=\frac {F_3(z)}{2z},\text{\,and\,}
F_1(z)=F_2(z)=\frac {4+3F_3(z)}{12z-9}.$$
In addition,
$ F_3(z)= {8z(4z-1)}/{(96z^3-72z^2-9)}.$
Therefore $$E(z^{-\tau})=f(z)=f_A(z)+f_B(z)+f_C(z)=\frac {16z^2-1}{3z(32z^3-24z^2-3)}.$$
Writing $z=1/\alpha$ yields that
$E(\alpha^\tau)=\frac {\alpha^2(\alpha^2-16)}{3(3\alpha^3+24\alpha-32)}$.
Taking $z=1$ gives
 ${f_A} = {f_B} = {1}/{{10}},{f_C} = {8}/{{10}},F_{1} = F_{2} = {44}/{15},F_{3} = {24}/{15}$, and hence
$E(\tau) = 1 + F_{1} + F_{2} + F_{3} = {127}/{15}$.
These results are all in agreement with that in  \cite{pozdnyakov2008occurrence} and \cite{gava2014stopping}.

Another way is to apply Corollary \ref{c2.3} and Corollary \ref{c2.4}. Because the first elements of $A,B,C$ are equal,
substituting
$${\tilde g_{AA}} = 17,{\tilde g_{BA}} = 1,{\tilde g_{CA}} = 1$$
$${\tilde g_{AB}} = 1,{\tilde g_{BB}} = 17,{\tilde g_{CB}} = 1$$
$${\tilde g_{AC}} = 1,{\tilde g_{BC}} = 1,{\tilde g_{CC}} = 3$$
into (\ref{eq1}) yields the following linear system:
$$\left\{ \begin{array}{l}
{f_A} + {f_B} + {f_C} = 1\\
- 16 \cdot {f_A} + 16 \cdot {f_B} = 0\\
- 16 \cdot {f_A}  + 2 \cdot {f_C} = 0
\end{array} \right.$$
Thus ${f_A} = {f_B} = {1}/{{10}}$ and $f_C = {8}/{{10}}$.
It is easy to see that the stationary distribution  is $\pi_1=\pi_2=\pi_3=1/3$.
Because the last elements of $A,B,C$ are all equal to $3$, by Corollary \ref{c2.4},
$$E(\tau|Z_1=3)=1+(f_A\cdot{\tilde g}_{AA}+ f_B\cdot{\tilde g}_{BA}+f_C\cdot{\tilde g}_{CA}-1)/\pi_{3}=29/5.$$
Clearly, $P(\tau_3=1)=\frac 13$ and $P(\tau_3=n)=\frac 23 \cdot(\frac 34)^{n-2}\cdot\frac 14$ for $n\ge 2$.
Therefore $$E(\tau)=E(\tau_3)-1+ E(\tau|Z_1=3)=127/15.$$
\end{example}

\begin{example}\label{exmp2}
Suppose that $\Delta {\rm{ = }}\left\{ {1,2} \right\},\mathcal C {\rm{ = }}\left\{ {A,B} \right\},A = 22,B = 121$ and
$$
P=\left (                %左括号
  \begin{array}{cc}
1/4& 3/4\\
3/4 & 1/4\\
  \end{array}
\right ).
$$
When will the distribution of $Z_\tau$ is the same as the initial distribution?

By calculating, we get
${ \tilde g_{AA}} = 5$, ${\tilde g_{BA} }= 0$,
${\tilde g_{AB}} = 0$ and ${\tilde g_{BB}} = {25}/{9}$.
The stationary distribution is $\pi_1=\pi_2=1/2$.
Using Corollary \ref{c2.4}, we have
$$\left\{ \begin{array}{l}
f_A+f_B=1\\
4\cdot f_A= \frac 12\cdot c\\
\frac {16}{9}\cdot f_B=\frac 12 \cdot c
\end{array} \right.$$
Hence $\mu_2={f_A} = {4}/{{13}},\mu_1={f_B} = {{9}}/{{13}}$ and $c = {{32}}/{{13}}$.
In addition,  $F_1=c\cdot \pi_1=16/13$, $F_2=c\cdot \pi_2=16/13$ and $E(\tau)=1+c=45/13$.
\end{example}

%The following is the appendix file.

%\attachfile{Appendix.pdf}

\end{document}